\documentclass[microtype]{gtpart}

\usepackage{latexsym,enumitem}
\usepackage{amssymb}
\usepackage[cp850]{inputenc}
\usepackage{epsfig}
\usepackage{psfrag}
\usepackage{amsthm}
\usepackage{amscd}
\usepackage{amsmath}
\usepackage{amsfonts}
\usepackage{graphics,caption}
\usepackage[all]{xy}
\usepackage{etoolbox}
\patchcmd{\quote}{\rightmargin}{\leftmargin 2em \rightmargin}{}{}
\usepackage[sort,numbers]{natbib} 
\usepackage{mathtools}
\usepackage{enumitem}

\usepackage{overpic}

\usepackage[english]{babel}
\usepackage{comment}
\usepackage{color}
\usepackage{hyperref}

\let\phi\varphi

\newcommand{\f}[2]{\frac{#1}{#2}} 

\let\epsilon\varepsilon

\newcommand{\be}{\begin{equation*}}
 \newcommand{\ee}{\end{equation*}}
\newcommand{\bpf}{\begin{dimo}}
\newcommand{\epf}{\end{dimo}}
\newcommand{\bdefi}{\begin{defin}}
\newcommand{\edefi}{\end{defin}}
\newcommand{\bthm}{\begin{thm}}
\newcommand{\ethm}{\end{thm}}
\newcommand{\blem}{\begin{lem}}
\newcommand{\elem}{\end{lem}}
\newcommand{\bcor}{\begin{cor}}
\newcommand{\ecor}{\end{cor}}

\newcommand{\bprop}{\begin{prop}}
\newcommand{\eprop}{\end{prop}}
\newcommand{\bese}{\begin{ese}}
\newcommand{\eese}{\end{ese}}
\newcommand{\brem}{\begin{rem}}
\newcommand{\erem}{\end{rem}}
\newcommand{\bpfc}{\begin{dimoclaim}}
\newcommand{\epfc}{\end{dimoclaim}}

\newcommand{\set}[1]{\left\{#1\right\}}					
	
\newenvironment{quot}
{
	\vspace{-0.2cm}
	\vspace{0.2cm}
}

\theoremstyle{definition}
\newtheorem{d1}{Definition}

\newenvironment{defin}
{
	\begin{quot}
		\begin{d1}
		}
		{\end{d1}
	\end{quot}

}

\theoremstyle{definition}
\newtheorem{r1}[d1]{Remark}

\newenvironment{rem}
{
	\begin{quot}
		\begin{r1}
		}
		{\end{r1}
	\end{quot}
}

\theoremstyle{definition}
\newtheorem{e1}[d1]{Exercise}

\theoremstyle{definition}
\newtheorem{ese1}[d1]{Example}

\newenvironment{ese}
{
	\begin{quot}
		\begin{ese1}
	}
	{	
		\end{ese1}
	\end{quot}
}

\theoremstyle{definition}

\theoremstyle{definition}
\newtheorem{f2}[d1]{Fact}

\theoremstyle{definition}

\theoremstyle{definition}

\theoremstyle{definition}
\newtheorem{t1}[d1]{Theorem}

\newenvironment{thm}
{
	\begin{quot}
		\begin{t1}}
		{\end{t1}
	\end{quot}
}

\theoremstyle{definition}
\newtheorem*{T1*}{Theorem}

\newenvironment{teor*}
{
	\begin{quot}
		\begin{T1*}}
		{\end{T1*}
	\end{quot}
}

\newenvironment{dimo}
{\begin{proof}[Proof]
	}
	{\end{proof}}

\newenvironment{dimoclaim}{\emph{Proof of Claim:}\;}{\hfill$\square$}

	\theoremstyle{definition}
	\newtheorem{l1}[d1]{Lemma}
	
	\newenvironment{lem}
	{
		\begin{quot}
			\begin{l1}}
			{\end{l1}
		\end{quot}
	}
	\theoremstyle{definition}
	\newtheorem{p1}[d1]{Proposition}
	
	\newenvironment{prop}
	{
		\begin{quot}
			\begin{p1}}
			{\end{p1}
		\end{quot}
	}
	
	\theoremstyle{definition}
	\newtheorem{c1}[d1]{Corollary}
	
	\newenvironment{cor}
	{
		\begin{quot}
			\begin{c1}}
			{\end{c1}
		\end{quot}
	}

		\renewenvironment{abstract}
	{\list{}{\rightmargin\leftmargin}%
		\item[\textbf{Abstract:}]\relax}
	{\endlist}

 \newtheorem*{Theorem*}{Theorem}
 \newtheorem*{Proposition*}{Proposition}
 \newtheorem*{Lemma*}{Lemma}	
  \newtheorem*{Corollary*}{Corollary}

\usepackage{multicol}

\usepackage{amsmath,amssymb,amsthm}
\usepackage{thmtools}
\usepackage{thm-restate}

\makeatletter
\DeclareFontFamily{OMX}{MnSymbolE}{}
\DeclareSymbolFont{MnLargeSymbols}{OMX}{MnSymbolE}{m}{n}
\SetSymbolFont{MnLargeSymbols}{bold}{OMX}{MnSymbolE}{b}{n}
\DeclareFontShape{OMX}{MnSymbolE}{m}{n}{
    <-6>  MnSymbolE5
   <6-7>  MnSymbolE6
   <7-8>  MnSymbolE7
   <8-9>  MnSymbolE8
   <9-10> MnSymbolE9
  <10-12> MnSymbolE10
  <12->   MnSymbolE12
}{}
\DeclareFontShape{OMX}{MnSymbolE}{b}{n}{
    <-6>  MnSymbolE-Bold5
   <6-7>  MnSymbolE-Bold6a
   <7-8>  MnSymbolE-Bold7
   <8-9>  MnSymbolE-Bold8
   <9-10> MnSymbolE-Bold9
  <10-12> MnSymbolE-Bold10
  <12->   MnSymbolE-Bold12
}{}

\let\llangle\@undefined
\let\rrangle\@undefined
\DeclareMathDelimiter{\llangle}{\mathopen}%
                     {MnLargeSymbols}{'164}{MnLargeSymbols}{'164}
\DeclareMathDelimiter{\rrangle}{\mathclose}%
                     {MnLargeSymbols}{'171}{MnLargeSymbols}{'171}
\makeatother

\renewcommand{\setminus}{\smallsetminus}

\begin{document}
\title[Fibered 3-manifolds with unique incompressible surfaces]{Fibered 3-manifolds with unique incompressible surfaces}
\author{Tommaso Cremaschi}
\address{University of Luxembourg, Department of Mathematics, Maison du Nombre, 6, avenue de la Fonte, L-4364 Esch-sur-Alzette, Grand Duchy of Luxembourg}
\email{tommaso.cremaschi@uni.lu}

\author{Andrew Yarmola}
\address{Department of Mathematics, PO Box 208283, 06520-8283, New Haven, CT, USA}
\email{andrew.yarmola@yale.edu}

\date{v1, \today}

\thanks{The first author acknowledges support from XXX.}

\maketitle
\begin{abstract}
We present a family $M_g$ of fibered hyperbolic 3-manifolds whose fibre $F_g$ is the unique connected incompressible surface and the genus $g \geq 2$ of $F_g$ can be arbitrary. This answers a question of Agol. \end{abstract}

\section{Introduction}
In \cite{AgolSmall}, Agol asks whether there exists ``fibred manifolds with fiber of arbitrarily high genus, such that the fiber is the only closed connected incompressible surface.'' In the present note, we answer this question in the affirmative by constructing closed (hyperbolic) fibered 3-manifolds with this property for every genus $g \geq 2$.

Our construction is based on the well known work of Hatcher and Thurston \cite{HT1985}, where the incompressible and $\partial$-incompressible surfaces in two-bridge knots are classified. We also make use of the work of Curtis-Franczak-Leiser-Manheimer \cite{CFLM2015} to understand when two such surfaces have the same boundary slope.

\bthm\label{main} For every $g \geq 2$, there exists a closed 3-manifold $M_g$ fibered over $\mathbb{S}^1$ where the fibre $F_g$ has genus $g$  and is the unique connected incompressible surface.
\ethm

Our family $\set{M_g}_{g \geq 2}$ is given by $0$-surgery on a family of fibered 2-bridge knots that have a unique essential surface with slope $0$, namely the fibre.

After writing this note, we were informed that two examples with $g = 2$ were posted on MathOverflow by Hatcher \cite{mathoverflow} in response to a question asking for manifolds where all incompressible surfaces are non-separating. This is true for our examples as well.

\section{Background}

A connected, incompressible, $\partial$-incompressible, and non-$\partial$-parallel properly embedded surface in $3$-manifold will be called \emph{essential}. In \cite{HT1985}, all essential surfaces in the exteriors of 2-bridge knots are classified. Further, it is shown that the only closed incompressible embedded surfaces are boundary parallel tori. 

For consistency, we will follow the notation and conventions of \cite{CFLM2015}, which differs slightly from \cite{HT1985}. A 2-bridge knot $K(\alpha, \beta)$ can be reconstructed from any signed continued fraction expansion:
\[
\f \beta\alpha=
r+\cfrac{1}{n_1+\cfrac{1}{n_2+\cfrac{1}{\ddots+\cfrac{1}{n_t}}}}= r + [n_1, n_2,\dotsc, n_t].\]
Since the $n_i$ are allowed to be negative, there are many such expansions. As shown in \cite{HT1985}, expansions with all $|n_i| \geq 2$ give rise to branched surfaces $\Sigma[n_1, \ldots, n_t]$ that carry all essential surfaces in $\mathbb{S}^3 \setminus \overset{\circ}{\mathcal{N}}(K(\alpha, \beta))$. Further, the authors of \cite{CFLM2015} demonstrate how all such expansions arise from the unique expansion $[m_1, \ldots, m_k]$ where $m_i > 0$ and $m_k \geq 2$. The correspondence is encoded in \emph{allowable sub-tuples} of indices $(i_1, \ldots, i_\ell)$ which are ordered sub-tuples of $(1, \ldots, k)$ satisfying certain conditions.

Essential surfaces carried by the same $\Sigma[n_1, \ldots, n_t]$ have the same boundary slope. In fact, this slope is always integral and can be computed from the corresponding allowable sub-tuple $(i_1, \ldots, i_\ell)$. The slope is given as \[m(n_1, \ldots, n_t) = \frac{1}{2}\left(\mathfrak c(i_1, \ldots, i_\ell) - \mathfrak c_0\right) \text{ where } \mathfrak c(i_1, \ldots, i_\ell) = \sum_{j = 1}^\ell (-1)^{i_j} m_{i_j}\]
and $\mathfrak c_0$ is the value of $\mathfrak c$ on an allowable sub-tuple corresponding to any branched surface carrying an orientable Seifert surface of $K(\alpha, \beta)$.

\section{Proof of Theorem \ref{main}}

\bpf Consider the 2-bridge knots $K_g = K(6g -1 , 2g)$ for $g \geq 2$. Then $[m_1, m_2, m_3] = [2,1, 2g-1]$ is the unique positive continued fraction with the last term $\geq 2$.

Following \cite{CFLM2015}, the allowable sub-tuples are $(1), (2),(3),$ and $(1,3)$. 
\begin{itemize}
\item $(1)$ corresponds to $\Sigma[-2,2, 2g-1]$ and $\mathfrak c(1) = -2$;
\item $(2)$  corresponds to $\Sigma[3, -2g ]$ and $\mathfrak c(2) = 1$;
\item  $(3)$ correspond to  $\displaystyle{\Sigma[2,2,\underbrace{-2,2}_{g - 1 \text{ times}}]}$  and  $\mathfrak c(3) = 1-2g$;
\item $(1,3)$ corresponds to $\displaystyle{\Sigma[-2,3,\underbrace{-2,2}_{g - 1\text{ times}}]}$ and $\mathfrak c(1,3) = -1-2g$;
\end{itemize} 
Notice that since $g \geq 2$, none of the values of $\mathfrak c$ are equal. Thus, by Corollary \cite[4.2]{CFLM2015}, all the slopes of essential surfaces carried by distinct branched surfaces are distinct.  Moreover, by Corollary of Proposition 1 in \cite{HT1985},  since $2g/(6g-1)$ can be expressed using only $\pm 2$, we see that $K_g$ is fibered and all orientable essential Seifert surfaces are isotopic to the fibre. The sub-tuple $(3)$ corresponds to this unique fibre. Combining these facts, we see that $\mathbb{S}^3 \setminus \overset{\circ}{\mathcal{N}}(K_g)$ has a unique essential surface with slope $0$. Doing $0$-filling gives an irreducible 3-manifold $M_g$ that is fibered and has a unique incompressible surface $F_g$. By \cite{GK1990}, the genus of $F_g$ is given by half the number of sign flips in the $\pm 2$ expansion, which is $g$. Note, one can read off the monodromy, as in \cite{GK1990}.

To see that $M_g$ is hyperbolic, observe that $M_g$ is irreducible and contains no essential tori. Further, since $H_1(M_g) = \mathbb Z$, we know by geometrization that $M_g$ is either hyperbolic or a Haken Seifert fibered space. However, $M_g$ contains an essential surface of genus $g \geq2$ and every such Haken Seifert fibered space must contain essential (vertical) tori. Thus, $M_g$ must be hyperbolic.
\epf

\bibliographystyle{plain}
\bibliography{mybib}

\end{document}